\documentclass[11pt]{article}
\usepackage[usenames]{color}
\usepackage{amssymb,amsmath,amsthm}

\numberwithin{equation}{section}

\newtheorem{thm}{Theorem}[section]

\newtheorem{lem}{Lemma}[section]

\theoremstyle{definition}
\newtheorem{defn}{Definition}[section]

\theoremstyle{remark}
\newtheorem{rem}{\bf Remark}[section]

\newcommand{\N}{\mbox{$\mathbb{N}$}}

\newcommand{\Z}{\mbox{$\mathbb{Z}$}}

\newcommand{\F}{\mbox{$\mathbb{F}$}}

\def \dv{\,|\,}
\def\lcm{\textrm{lcm}}

\begin{document}

\vskip 1 true cm

\centerline{\bf  On the counting function of sets with even
partition functions} \centerline{by}

\centerline {\bf F. Ben Sa\"id }\centerline{Universit\'e de
Monastir} \centerline{Facult\'e des Sciences de Monastir}
\centerline{ Avenue de l'environnement, 5000 Monastir, Tunisie.}
\centerline{Fethi.BenSaid@fsm.rnu.tn} \centerline {and} \centerline
{\bf J.-L. Nicolas } \centerline{Universit\'e de Lyon, Universit\'e
Lyon 1, CNRS} \centerline{Institut Camile Jordan, Math\'ematiques}
\centerline{ Batiment Doyen Jean Braconnier}
\centerline{Universit\'e Claude Bernard} \centerline{ 21 Avenue
Claude Bernard, F-69622 Villeurbanne cedex, France}
\centerline{jlnicola@in2p3.fr} \vskip 1 true cm

\hskip 7 true cm
\begin{minipage}[t]{7cm}
\noindent {\bf To K\'alm\'an Gy\H ory, Attila Peth\H o,  J\'anos
Pintz and Andr\'as S\'ark\"ozy for their nice works in number
theory.}
\end{minipage}
\vskip 1 true cm
\begin{abstract}
Let $q$ be an odd positive integer and $P\in \F_2[z]$ be of order
$q$ and such that $P(0)=1$. We denote by ~${\cal A}={\cal A}(P)$ the
unique set of positive integers satisfying
$\sum_{n=0}^{\infty}p({\cal A},n)z^n\equiv P(z) ~(\bmod ~2)$, where
$p({\cal A},n)$ is the number of partitions of $n$ with parts in
${\cal A}$. In \cite{F lalouar}, it is proved that if $A(P,x)$ is
the counting function of the set ${\cal A}(P)$ then $A(P,x)\ll
x(\log x)^{-r/\varphi (q)}$, where $r$ is the order of $2$ modulo
$q$ and $\varphi$ is the Euler's function. In this paper, we improve
on the constant $c=c(q)$ for which $A(P,x)\ll x(\log x)^{-c}$.
\end{abstract}

key words: Sets with even partition functions, bad and semi-bad primes, order of a polynomial, Selberg-Delange formula.\\

\section{Introduction.}
Let $\N$ be the set of positive integers and ${\cal A}=\{ a_1,a_2,
...\}$ be a subset of $\N$. For $n\in \N$, we denote by $p({\cal
A},n)$ the number of partitions of $n$ with parts in ${\cal A}$,
i.e. the number of solutions of the equation
$$a_1x_1+a_2x_2+ ...=n,$$
in non-negative integers $x_1,x_2, ...$.  We set $p({\cal A},0)=1$.\\

Let $\F_2$ be the field with two elements and $f=1+\epsilon_1
z+...+\epsilon_Nz^N+\cdots \in \F_2[[z]]$. Nicolas et al proved (see
\cite {NRS}, \cite {Ben2} and \cite{Lah}) that there is a unique
subset ${\cal A}={\cal A}(f)$ of $\N$ such that
\begin{equation}\label{cong}
\sum_{n=0}^\infty p({\cal A},n)z^n\equiv f(z)~(\bmod~2).
\end{equation}
When $f$ is a rational fraction, it has been shown in \cite{Lah}
that there is a polynomial $U$ such that ${\cal A}(f)$ can be easily
determined from ${\cal A}(U)$. When $f$ is a general power series,
nothing  about the behaviour of ${\cal A}(f)$ is known. From now on,
we shall restrict ourselves to the case $f=P$, where
$$P=1+\epsilon_1 z+...+\epsilon_Nz^N\in \F_2[z]$$
is a polynomial of degree $N\geq 1$.\\

Let $A(P,x)$ be the counting function of the set ${\cal A}(P)$, i.e.
\begin{equation}\label{count}
A(P,x)=\mid \{n: 1\leq n\leq x, n\in {\cal A}(P) \}\mid.
\end{equation}
In \cite{Dai}, it is proved that
\begin{equation}\label{Dai}
A(P,x)\geq \frac {\log x}{\log 2}-\frac {\log (N+1)}{\log 2}.
\end{equation}
More attention was paid on upper bounds for $A(P,x)$. In
\cite[Theorem
3]{F lalouar}, it was observed that when $P$ is a product of cyclotomic polynomials, the set ${\cal A}(P)$ is a union of geometric progressions of quotient $2$ and so $A(P,x)={\cal O}(\log x)$.\\

Let the decomposition of $P$ into irreducible factors over $\F_2[z]$
be
$$P=P_1^{\alpha_1}P_2^{\alpha_2}\cdots P_l^{\alpha_l}.$$
We denote by $\beta_i$, $1\leq i\leq l$, the order of $P_i(z)$, that
is the smallest positive integer such that $P_i(z)$ divides
$1+z^{\beta_i}$ in $\F_2[z]$; it is known that $\beta_i$ is odd (cf.
\cite{LN}). We set
\begin{equation}\label{q^P}
q=q(P)=\lcm (\beta_1, \beta_2, ..., \beta_l).
\end{equation}
If $q=1$ then $P(z)=1+z$ and ${\cal A}(P)=\{2^k,~k\geq 0\}$, so that
$A(P,x)={\cal O}(\log x)$. We may suppose that $q\geq 3$. Now, let
\begin{equation}\label{sigma}
\sigma ({\cal A}, n)=\sum_{d\dv n,~d\in {\cal A}}d=\sum_{d\dv
n}d\chi({\cal A},d),
\end{equation}
where $\chi({\cal A},.)$ is the characteristic function of the set
${\cal A}$,
$$\chi({\cal A},d)=\left\{
              \begin{array}{l}
               1~\text {if}~d\in {\cal A}\\
                0~\text {otherwise}.
              \end{array}
            \right.$$
In \cite{Ben3} (see also \cite{Ben1} and \cite{NFA}), it is proved
that for all $k\geq 0$, $q$ is a period of the sequence $(\sigma
({\cal A}, 2^kn)~\bmod~2^{k+1})_{n\geq 1}$, i.e.
\begin{equation}\label{period}
n_1\equiv n_2 ~(\bmod~q)\Rightarrow \sigma ({\cal A}, 2^kn_1)~\equiv
\sigma ({\cal A}, 2^kn_2)~(\bmod~2^{k+1})
\end{equation}
and $q$ is the smallest integer such that (\ref{period}) holds for
all $k's$. Moreover, if $n_1$ and $n_2$ satisfy  $n_2~\equiv
2^an_1~(\bmod~q)$ for some $a\geq 0$, then
\begin{equation}\label{periodbis}
\sigma ({\cal A}, 2^kn_2)~\equiv \sigma ({\cal A},
2^kn_1)~(\bmod~2^{k+1}).
\end{equation}
If $m$ is odd and $k\geq 0$, let
\begin{equation}\label{S}
S_{{\cal A}}(m,k)=\chi({\cal A},m)+2\chi({\cal
A},2m)+\ldots+2^k\chi({\cal A},2^km).
\end{equation}
It follows that for $n=2^km$, one has
\begin{equation}\label{sigma2}
\sigma ({\cal A}, n)=\sigma ({\cal A}, 2^km)=\sum_{d\dv m}dS_{{\cal
A}}(d,k),
\end{equation}
which, by M\"obius inversion formula, gives
\begin{equation}\label{inversion}
mS_{{\cal A}}(m,k)=\sum_{d\dv m}\mu (d)\sigma ({\cal A}, \frac
{n}{d})= \sum_{d\dv \overline{m}}\mu (d)\sigma ({\cal A}, \frac
{n}{d}),
\end{equation}
where $\mu$ is the M\"obius's function and $\overline{m}=\prod_{p\dv m}p$ is the radical of $m$, with $\overline{1}=1$.\\

In \cite{Ben4} and \cite{FJZE}, precise descriptions of the sets
${\cal A}(1+z+z^3)$ and ${\cal A}(1+z+z^3+z^4+z^5)$ are given and
asymptotics to the related counting functions are obtained,
\begin{equation}\label{123}
A(1+z+z^3,x)\sim c_1\frac {x}{(\log x)^{\frac
{3}{4}}},~~x\rightarrow \infty,
\end{equation}
\begin{equation}\label{12345}
A(1+z+z^3+z^4+z^5,x)\sim c_2\frac {x}{(\log x)^{\frac
{1}{4}}},~~x\rightarrow \infty,
\end{equation}
where $c_1=0.937..., c_2=1.496...$. In \cite{NF2}, the sets ${\cal A}(P)$ 
are considered when $P$ is irreducible of prime order $q$
and such that the order of $2$ in $(\Z/q\Z)^*$ is $\frac {q-1}{2}$. 
This situation is similar to that of ${\cal A}(1+z+z^3)$, and 
formula (\ref{123}) can
be extended to $A(P,x)\sim c' x(\log x)^{-3/4},~x\rightarrow \infty $, 
for some constant $c'$ depending on $P$.\\

Let $P=QR$ be the product of two coprime polynomials in $\F_2[z]$.
In \cite{Ben2}, the following is given
\begin{equation}\label{A'A"}
A(P,x)\leq A(Q,x)+A(R,x)
\end{equation}
and
\begin{equation}\label{A"}
\mid A(P,x)-A(R,x)\mid \leq \sum_{0\leq i\leq \frac{\log x}{\log
2}}A(Q,\frac {x}{2^i}) .
\end{equation}
As an application of (\ref{A"}), choosing $Q=1+z+z^3$,
$R=1+z+z^3+z^4+z^5$ and $P=QR$, we get from (\ref{123})-(\ref{A"}),
$$A(P,x)\sim A(R,x)\sim c_2 x(\log x)^{-1/4}, ~x\rightarrow \infty .$$

In \cite{F lalouar}, a claim of Nicolas and S\'ark\"ozy
\cite{NS2001}, that some polynomials with $A(P,x)\asymp x$ may
exist, was disapproved. More precisely, the following was obtained
\begin{thm}\label{thm} Let $P\in \F_2[z]$ be such that $P(0)=1$,
  ${\cal A}={\cal A}(P)$
be the unique set obtained from (\ref{cong}) and $q$ be the odd number
defined by (\ref{q^P}). Let $r$ be the order of $2$ modulo $q$, that
is  the smallest positive integer such that $2^r\equiv 1~(\bmod~q)$.
We shall say that a prime $p\not =2$ is a bad prime if
\begin{equation}\label{badprime}
\exists~~ i, ~~0\leq i\leq r-1~~\text {and}~~p\equiv 2^i ~(\bmod~q).
\end{equation}
(i) If $p$ is a bad prime, we have $\gcd(p,n)=1$ for all $n\in {\cal A}$.\\
(ii) There exists an absolute constant $c_3$ such that for all
$x>1$,
\begin{equation}\label{upperbound}
A(P,x)\leq 7(c_3)^r \frac {x}{(\log x)^{\frac {r}{\varphi (q)}}},
\end{equation}
where $\varphi $ is Euler's function.
\end{thm}
\section{The sets of bad and semi-bad primes.}
Let $q$ be an odd integer $\geq 3$ and $r$ be the order of $2$
modulo $q$. Let us call "bad classes" the elements of
\begin{equation}\label{calE(q)}
{\cal E}(q)=\{ 1, 2, ..., 2^{r-1}\}\subset (\Z/q\Z)^*.
\end{equation}
From (\ref{badprime}), we know that an odd prime $p$ is bad if
$p~\bmod~q$ belongs to ${\cal E}(q)$. The set of bad primes will be
denoted by ${\cal B}$. The fact that no element of ${\cal A}(P)$ is
divisible by a bad prime (cf. Theorem \ref{thm} (i)) has given (cf.
\cite{F lalouar}) the upper bound (\ref{upperbound}). Two other
sets of primes will be used  to improve (\ref{upperbound}) cf. Theorem \ref{thimproved} below.\\
\begin{rem} $2$ is not a bad prime although it is a bad class.
\end{rem}
\begin{defn} A class of $(\Z/q\Z)^*$ is said semi-bad if it does not
  belong to ${\cal E}(q)$ and its square does. A prime $p$ is called
semi-bad if its class modulo $q$ is semi-bad. We denote by ${\cal
E}'(q)$ the set of semi-bad classes, so that
\begin{equation}\nonumber
p \;\text{ semi-bad} \quad \Longleftrightarrow \quad p~\bmod~q\in
{\cal E}'(q).
\end{equation}
We denote by $| {\cal E}'(q)|$ the number of elements of ${\cal
E}'(q)$.
\end{defn}

\begin{lem}\label{1}
Let $q$ be an odd integer $\geq 3$, $r$ be the order of $2$
modulo $q$ and
$$q_2=\left\{
              \begin{array}{l}
               1~\text {\rm  if $2$ is a square modulo q}\\
                0~\text {\rm  if not}.
                \end{array}
                 \right.$$
The number $|{\cal E}'(q)|$ of semi-bad classes modulo $q$ is given by
\begin{eqnarray}\label{Sq}
|{\cal E}'(q)| & = & 2^{\omega (q)}
\left(\left\lfloor  \frac {r+1}{2}\right\rfloor+q_2\left\lfloor \frac
    {r}{2}\right\rfloor \right)-r\notag\\
& = & \begin{cases}
r(2^{\omega(q)-1}-1) & \text{ if $r$ is even and $q_2=0$}\\
r(2^{\omega(q)}-1) & \text{otherwise},\\
\end{cases}
\end{eqnarray}
where $\omega (q)$ is the number of distinct prime factors of $q$ and $\left\lfloor x\right\rfloor $ is the floor of $x$.
\end{lem}
\begin{proof} We have to count the number of solutions of the $r$ congruences
\begin{equation}\nonumber
E_i:\ \ x^2\equiv 2^i~(\bmod~q),\ \ 0\leq i\leq r-1,
\end{equation}
which do not belong to ${\cal E}(q)$. The number of solutions of
$E_0$ is $2^{\omega (q)}$. The contribution of $E_i$ when $i$ is
even is equal to that of $E_0$ by the change
of variables $x=2^{i/2}\xi $, so that the total number of solutions, in $(\Z/q\Z)^*$, of the $E_i's$ for $i$ even is equal to $\left\lfloor  \frac {r+1}{2}\right\rfloor 2^{\omega (q)} $.\\

The number of odd $i's$, $0\leq i\leq r-1$, is equal to $\left\lfloor \frac {r}{2}\right\rfloor $. The contribution of all the $E_i's$ for these $i's$ are equal
and vanish if $q_2=0$. When $q_2=1$, $E_1$
has $2^{\omega (q)}$ solutions in $(\Z/q\Z)^*$. Hence the total number of solutions, in $(\Z/q\Z)^*$, of the $E_i's$ for $i$ odd is equal
to $q_2\left\lfloor  \frac {r}{2}\right\rfloor 2^{\omega (q)} $.\\

Now, we have to remove those solutions which are in ${\cal E}(q)$.
But any element $2^i$, $0\leq i\leq r-1$, from ${\cal E}(q)$ is a
solution of the congruence $x^2\equiv 2^j~(\bmod~q)$, where
$j=2i~\bmod~r$. Hence
\begin{equation}\nonumber
{|\cal E}'(q)|=2^{\omega (q)}\left(\left\lfloor  \frac {r+1}{2}\right\rfloor+q_2\left\lfloor \frac {r}{2}\right\rfloor \right)-r.
\end{equation}
The second formula in \eqref{Sq} follows by noting that $q_2=1$ when
$r$ is odd.
\end{proof}

\begin{defn} A set of semi-bad classes
is called a coherent set if
it is not empty and if the product of
any two  of its elements is a bad class.
\end{defn}

\begin{lem}\label{lemsb}
Let $b$ be a semi-bad class; then
$${\cal C}_b=\{b,2b,\ldots,2^{r-1}b\}$$
is a coherent set. There are no coherent sets with more than $r$
elements.
\end{lem}

\begin{proof}
First, we observe that, for $0\leq u \leq r-1$, $2^ub$ is semi-bad and,
for $0\leq u < v \leq r-1$, $(2^ub)(2^vb)$ is bad so that ${\cal C}_b$
is coherent.

Further, let ${\cal F}$ be a set of semi-bad classes with more than
$r$ elements; there exists in ${\cal F}$ two semi-bad classes $a$
and $b$ such that $a\notin {\cal C}_b$. Let us prove that $ab$ is
not bad. Indeed, if $ab\equiv 2^u ~(\bmod~ q)$ for some $u$, we
would have $a\equiv 2^ub^{-1} ~(\bmod~q)$. But, as $b$ is semi-bad,
$b^2$ is bad, i.e. $b^2\equiv 2^v ~(\bmod~q)$ for some $v$, which
would imply $b\equiv 2^vb^{-1}~(\bmod~q)$, $b^{-1}\equiv
b2^{-v}~(\bmod~q)$, $a\equiv 2^{u-v}b~(\bmod~q)$ and $a\in {\cal
  C}_b$, a contradiction. Therefore, ${\cal F}$ is not coherent.
\end{proof}

\begin{lem}\label{lemE'}
If $\omega(q)=1$ and $\varphi(q)/r$ is odd, then ${\cal
  E}'(q)=\emptyset$; while if $\varphi(q)/r$ is even, the set of semi-bad
classes ${\cal E}'(q)$ is a coherent set of $r$ elements.

If $\omega(q)\geq 2$, then ${\cal E}'(q) \neq \emptyset$ and there
exists
  a coherent set ${\cal C}$ with $|{\cal C}|=r$.
\end{lem}

\begin{proof}
If $\omega(q)=1$, $q$ is a power of a prime number and the group
$(\Z/q\Z)^*$ is cyclic. Let $g$ be some generator and $d$ be the
smallest positive integer such that $g^d \in {\cal E}(q)$, where
${\cal E}(q)$ is given by (\ref{calE(q)}). We have $d=\varphi(q)/r$,
since $d$ is the order of the group $(\Z/q\Z)^*/_{{\cal E}(q)}$. The
discrete logarithms of the bad classes are $0,d,2d,\cdots ,(r-1)d$.
The set ${\cal E}'(q) \cup {\cal E}(q)$ is equal to the union of the
 solutions of the congruences
\begin{equation}\label{congq}
x^2\equiv g^{ad}\pmod{q}
\end{equation}
for $0\leq a \leq r-1$. By the change of variable $x=g^t$,
\eqref{congq} is equivalent to
\begin{equation}\label{congphiq}
2t\equiv ad \pmod{\varphi(q)}.
\end{equation}

Let us assume first that  $d$ is odd so that $r$ is even.  If $a$ is
odd, the congruence \eqref{congphiq} has no solution while, if $a$ is
even, say $a=2b$, the solutions of \eqref{congphiq} are
$t\equiv bd \pmod{\varphi(q)/2}$ i.e.
$$t\equiv bd \pmod{\varphi(q)} \quad \text{ or } \quad   t\equiv bd
+(r/2)d \pmod{\varphi(q)},$$ which implies
$${\cal E}'(q) \cup {\cal E}(q)=\{g^0,g^d,\ldots,g^{(r-1)d}\}={\cal E}(q)$$
and ${\cal E}'(q)=\emptyset$.

Let us assume now that $d$ is even.  The congruence \eqref{congphiq}
is equivalent to
$$t\equiv ad/2 \pmod{\varphi(q)/2}$$
which implies ${\cal E}'(q) \cup {\cal E}(q)=\{g^{\alpha d/2}, 0\leq
\alpha \leq 2r-1\}$ yielding
$${\cal E}'(q)=\{g^{\frac {d}{2}},g^{3\frac
  {d}{2}},\cdots ,g^{(2r-1)\frac {d}{2}}\}={\cal C}_b$$
(with $b=(g^{\frac {d}{2}})$), which is coherent by Lemma
\ref{lemsb}.\\

If $\omega(q)\geq 2$, then, by Lemma \ref{1},  ${\cal E}'(q) \neq
  \emptyset$. Let $b\in {\cal E}'(q)$; by Lemma \ref{lemsb}, the set
  ${\cal C}_b$ is a coherent set of $r$ elements.
\end{proof}

Let us set
\begin{equation}\label{c(q)}
 c(q)=
\begin{cases}
\frac {3}{2}  & \mbox{ if ${\cal E}'(q)\neq \emptyset$ }\\
1 & \mbox{ if ${\cal E}'(q)= \emptyset$. }\\
\end{cases}
\end{equation}
We shall prove
\begin{thm}\label{thimproved}
Let $P\in \F_2[z]$ with $P(0)=1$, $q$ be the odd integer defined by
(\ref{q^P}) and $r$ be the order of $2$ modulo $q$. We denote by ${\cal
A}(P)$ the set obtained from (\ref{cong}) and by $A(P,x)$ its
counting function. When $x$ tends to infinity, we have
\begin{equation}\label{rodd}
A(P,x)\ll_q \frac {x}{(\log x)^{c(q)\frac {r}{\varphi (q)}}},
\end{equation}
where $c(q)$ is given by (\ref{c(q)}).
\end{thm}
When $P$ is irreducible, $q$ is prime and
$r=\frac {q-1}{2}$, the upper bound \eqref{rodd} is best possible;
indeed  in this case, from \cite{NF2}, we have
$A(P,x)\asymp \frac {x}{(\log x)^{3/4}}$.
As $\varphi(q)/r=2$, Lemma \ref{lemE'} implies
${\cal E}'(q)\neq \emptyset$ so that $c=3/2$ and in \eqref{rodd}, the
exponent of $\log x$ is $3/4$.
Moreover, formula
(\ref{12345}) gives the optimality of (\ref{rodd}) for some prime
($q=31$) satisfying $r=\frac{q-1}{6}$.

\begin{thm}\label{thimproved2}
Let $P\in \F_2[z]$ be such that $P(0)=1$ and $P=P_1P_2\cdots P_j$,
where the $P_i's$ are irreducible polynomials in $\F_2[z]$. For
$1\leq i\leq j$, we denote by $q_i$ the order of $P_i$, by $r_i$ the
order of $2$ modulo $q_i$ and we set $c=\min_{1\leq i\leq j}c(q_i)r_i/\varphi(q_i)$,
where $c(q_i)$ is given by (\ref{c(q)}). When $x$ tends to infinity, we have
\begin{equation}\label{Min}
                  A(P,x)\ll\frac {x}{(\log x)^{c}} \cdot
                  \end{equation}
where the  symbol $\ll$ depends on the $q_i's$, $1\leq i\leq j$.
\end{thm}

Let ${\cal C}$ be a coherent set of semi-bad classes modulo $q$. Let
us associate to ${\cal C}$ the set of primes ${\cal S}$ defined by
\begin{equation}\label{CS}
p \in {\cal S}\quad \Longleftrightarrow\quad p~\bmod~ q \in {\cal C}.
\end{equation}
We define  $\omega_{\cal S}$ as the additive arithmetic function
\begin{equation}\label{wS}
\omega_{\cal S}(n)=\sum_{p\dv n,~ p\in {\cal S}}1.
\end{equation}

\begin{lem}\label{2}
Let $m$ be an odd positive integer, not divisible by
any bad prime. If $\omega_{\cal S}(m)= k+2\geq 2$ then $2^hm\not \in
{\cal A}(P)$ for all $h$, $0\leq h\leq k$. In other words, if
$2^hm \in {\cal A}(P)$, then $h \geq \omega_{\cal S}(m)-1$ holds.
\end{lem}

\begin{proof}
Let us write $\overline{m}=m'm"$, with $m'=\prod_{p\dv
\overline{m},~ p\in {\cal S}}p$ and $m"=\prod_{p\dv \overline{m},~
p\not \in {\cal S}}p$. From (\ref{inversion}), if $n=2^km$ then
\begin{equation}\label{inversion2}
mS_{{\cal A}}(m,k)=\sum_{d\dv \overline{m}}\mu (d)\sigma ({\cal A},
\frac {n}{d})=\sum_{d'\dv m'}\sum_{d"\dv m"}\mu (d')\mu (d")\sigma
({\cal A}, \frac {n}{d'd"}).
\end{equation}
Let us write $d'=p_{i_1}\cdots p_{i_j}$ and take some $p_{\cal S}$
from ${\cal S}$. If $j$ is even then $\mu (d')=1$ and, from
the definition of a coherent set, $d'\equiv
2^{t}~(\bmod~q)$ for some $t$ (depending on $d'$), $0\leq t\leq
r-1$. Whereas, if $j$ is odd then $\mu (d')=-1$ and $d'\equiv
2^{t'}p_{\cal S}^{-1}~(\bmod~q)$ for some $t'$ (depending on $d'$),
$0\leq t'\leq r-1$. From (\ref{periodbis}), we obtain

\begin{equation}\label{invers1}
\mu (d')\sigma ({\cal A}, \frac {n}{d'd"})\equiv \sigma ({\cal A},
\frac {n}{d"})~ (\bmod ~ 2^{k+1})~ \text {if $j$ is even},
\end{equation}
\begin{equation}\label{invers2}
\mu (d')\sigma ({\cal A}, \frac {n}{d'd"})\equiv -\sigma ({\cal A},
\frac {np_{\cal S}}{d"})~ (\bmod ~ 2^{k+1})~ \text {if $j$ is odd}.
\end{equation}
Since $\alpha =\omega_{\cal S}(\overline{m})=k+2 > 0$, the number of
$d'$ with odd $j$ is equal to that with even $j$ and is given by
\begin{equation}\nonumber
1+\binom {\alpha}{2}+\binom {\alpha}{4}+\cdots =\binom
{\alpha}{1}+\binom {\alpha}{3}+\cdots =2^{\alpha -1}.
\end{equation}
From (\ref{inversion2}), we obtain
\begin{equation}\label{inversion3}
mS_{{\cal A}}(m,k)\equiv2^{\alpha -1}\sum_{d"\dv m"}\mu (d")\left(
\sigma ({\cal A}, \frac {n}{d"})-\sigma ({\cal A}, \frac {np_{\cal
S}}{d"})\right)~ (\bmod ~ 2^{k+1}),
\end{equation}
which, as $\alpha =\omega_{\cal S}(m)= k+2$, gives $S_{{\cal
A}}(m,k)\equiv 0~ (\bmod ~ 2^{k+1})$, so that from (\ref{S}),
\begin{equation}\label{notin}
\chi({\cal A},m)=\chi({\cal A},2m)= \cdots =\chi({\cal A},2^km)=0.
\end{equation}
\end{proof}
Let us assume that ${\cal E}'(q)\neq \emptyset$ so that there exists
a coherent set  ${\cal C}$ with $r$ semi-bad classes modulo $q$; we
associate to ${\cal C}$ the set of primes ${\cal S}$ defined by
\eqref{CS} and we denote by ${\cal Q}={\cal Q}(q)$ and ${\cal
N}={\cal N}(q)$ the sets
\begin{equation}\nonumber
{\cal Q}=\{ p\text{ prime},\ p\dv q \}\ \ \text {and}\ \ {\cal N}=\{
p\text{ prime},\ p\not\in {\cal B}\cup {\cal S}\ \text {and}\ \gcd
(p,2q)=1 \},
\end{equation}
so that the whole set of primes is equal to ${\cal B}\cup {\cal
S}\cup {\cal N}\cup {\cal Q}\cup \{2\}$. For $n\geq 1$, let us
define the multiplicative arithmetic function
$$\delta(n)=
   \begin{cases}
         1 & \text{ if } ~p\dv n\Rightarrow p\not \in {\cal B}~(i.e. ~p\in {\cal S}\cup {\cal N}\cup {\cal Q}\cup \{2\}) \\
         0 & \text{ otherwise}.
    \end{cases}
$$
and for $x> 1$,
  \begin{equation}\label{Vx}
  V(x)= V_q(x)=\sum_{n\geq 1,~n2^{\omega_{\cal S}(n)}\leq x}\delta(n).
\end{equation}
\begin{lem}\label{3} Under the above notation, we have
\begin{equation}\label{Vxasym}
  V(x)=V_q(x)={\cal O}_q\left( \frac {x}{(\log x)^{c(q)\frac {r}{\varphi (q)}}}
\right),
\end{equation}
where $c(q)$ is given by (\ref{c(q)}).
\end{lem}
\begin{proof}
To prove (\ref{Vxasym}), one should consider, for complex $s$ with
${\cal R}(s)>1$, the series
\begin{equation}\label{F(s)}
F(s)=\sum_{n\geq 1}\frac {\delta (n)}{(n2^{\omega_{\cal S}(n)})^s}\cdot
\end{equation}
This  Dirichet series has an Euler's product given by
\begin{equation}\label{F(s)Euler}
F(s)=\prod_{p\in {\cal N}\cup {\cal Q}\cup \{ 2\}}\left( 1-\frac
{1}{p^s}\right)^{-1}~~\prod_{p\in {\cal S}}\left( 1+\frac
{1}{2^s(p^s-1)}\right),
\end{equation}
which can be written as
\begin{equation}\label{F(s)rewritten}
F(s)=H(s)\prod_{p\in {\cal N}}\left( 1-\frac
{1}{p^s}\right)^{-1}\prod_{p\in {\cal S}}\left( 1-\frac
{1}{p^s}\right)^{-\frac {1}{2^s}},
\end{equation}
where
\begin{equation}\label{G(s)}
H(s)=\prod_{p\in {\cal Q}\cup \{2\}}\left( 1-\frac
{1}{p^s}\right)^{-1} \prod_{p\in {\cal S}}\left( 1+\frac
{1}{2^s(p^s-1)}\right)\left( 1-\frac {1}{p^s}\right)^{\frac
{1}{2^s}}.
\end{equation}
By applying Selberg-Delange's formula (cf. \cite{Ben5}, Th\'eor\`eme
$1$  and \cite{FJZE}, Lemma $4.5$), we obtain some constant $c_4$
such that
\begin{equation}\label{V(x)asymp}
V(x)=c_4 \frac {x}{(\log x)^{c(q)\frac {r}{\varphi (q)}}}+{\cal
O}_q\left( \frac {x\log \log x}{\log x}\right).
\end{equation}
The constant $c_4$ is somewhat complicated, it is given by
\begin{equation}\label{c5}
c_4=\frac {C H(1)}{\Gamma (1-c(q)\frac {r}{\varphi (q)})},
\end{equation}
where $\Gamma $ is the gamma function,
\begin{equation}\label{G(1)}
H(1)=\frac {2q}{\varphi (q)}\prod_{p\in {\cal S}}\left(1+\frac
{1}{2(p-1)}\right)\left(1-\frac {1}{p}\right)^{\frac {1}{2}}
\end{equation}
and
\begin{equation}\nonumber
C=\prod_{p\in {\cal N}}\left( 1-\frac {1}{p}\right)^{-1}\prod_{p\in
{\cal S}}\left( 1-\frac {1}{p}\right)^{\frac {-1}{2}}\prod_{p}\left(
1-\frac {1}{p}\right)^{1-c(q)\frac {r}{\varphi (q)}},
\end{equation}
where in the third product, $p$ runs over all primes.
\end{proof}

\section{Proof of the results.}
{\bf Proof of Theorem \ref{thimproved}.} If $r=\varphi (q)$
then $2$ is a generator of $(\Z/q\Z)^*$, all primes
are bad but $2$ and the prime factors of $q$; hence
by Theorem $2$ of \cite{F lalouar},
$A(P,x)={\cal O}\left((\log x)^\kappa \right)$ for some
constant $\kappa$, so that we may remove the case $r=\varphi (q)$.

If ${\cal E}'(q)=\emptyset$, from \eqref{c(q)}, $c=1$
holds and \eqref{rodd} follows from \eqref{upperbound}.

We now assume ${\cal E}'(q)\neq \emptyset$, so that,
from Lemma \ref{lemsb}, there exists a coherent set ${\cal C}$
satisfying $|{\cal C}|=r$. We define the set of primes ${\cal S}$ by \eqref{CS}.
Let us write $V(x)$ defined in (\ref{Vx}) as
\begin{equation}\label{V'V"}
V(x)=V'(x)+V"(x),
\end{equation}
with
\begin{equation}\nonumber
V'(x)=\sum_{n\geq 1,\ n2^{\omega_{\cal S}(n)}\leq x,\ \omega_{\cal
S}(n)=0}\delta (n) \ \ \text {and}\ \ V"(x)=\sum_{n\geq 1,\
n2^{\omega_{\cal S}(n)}\leq x,\ \omega_{\cal S}(n)\geq 1}\delta (n).
\end{equation}
Similarly, we write $A(P,x)=\sum_{a\in {\cal A}(P),\ a\leq
x}1=A'+A"$, with
\begin{equation}\nonumber
A'=\sum_{a\in {\cal A}(P),\ a\leq x,\ \omega_{\cal S}(a)=0}1\ \
\text {and}\ \ A"=\sum_{a\in {\cal A}(P),\ a\leq x,\ \omega_{\cal
S}(a)\geq 1}1.
\end{equation}
An element $a$ of ${\cal A}(P)$ counted in $A'$ is free of bad and
semi-bad primes, so that
\begin{equation}\label{A'V'}
A'\leq V'(x)\leq V'(2x).
\end{equation}
By Lemma \ref{2}, an element $a$ of ${\cal A}(P)$ counted in $A"$ is
of the form $n2^{\omega_{\cal S}(n)-1}$ with
$\omega_{\cal S}(n) = \omega_{\cal S}(a) \geq 1$; hence
\begin{equation}\label{A"V"}
A"\leq V"(2x).
\end{equation}
Therefore, from (\ref{V'V"})-(\ref {A"V"}), we get
\begin{equation}\nonumber
A(P,x)=A'+A"\leq V'(2x)+V"(2x)=V(2x)
\end{equation}
and (\ref{rodd}) follows from Lemma \ref{3}. \hskip 9.5 true cm
$\Box $

$\ $\\
{\bf Proof of Theorem \ref{thimproved2}.} Just use Theorem \ref{thimproved} and (\ref {A'A"}).\hskip 4.7 true cm $\Box $\\

{\bf Acknowledgements.} It is a great pleasure for us to thank Andr\'as
S\'ark\"ozy for initiating the study of the sets ${\cal A}(P)$ and
for all the mathematics that we have learnt from him.
\def\refname{References}

\end{document}